\documentclass[11pt,a4paper]{amsart}

\usepackage{a4,amsmath,amsfonts,amssymb, amsthm,
latexsym,enumerate,xy,graphics}

\xyoption{all}
\CompileMatrices

\newcommand{\be}{\begin{equation}}
\newcommand{\ee}{\end{equation}}

\def\C{\ensuremath{\mathbb{C}}}

\def\P{\ensuremath{\mathbb{P}}}
\def\Q{\ensuremath{\mathbb{Q}}}

\def\Z{\ensuremath{\mathbb{Z}}}

\def\dim{\mathop{\mathrm{dim}}\nolimits}

\def\RHom{\mathop{\mathbf{R}\mathrm{Hom}}\nolimits}
\def\id{\mathop{\mathrm{id}}\nolimits}

\def\Lie{\mathop{\mathrm{Lie}}\nolimits}

\def\rk{\mathop{\mathrm{rk}}}

\def\Spec{\mathop{\mathrm{Spec}}}

\def\dbyd#1{\ensuremath{\frac{\partial}{\partial #1}}}

\newenvironment{Prf}{\textit{Proof.}\/}{\hfill$\Box$}

\newtheorem{Thm}{Theorem}[subsection]
\newtheorem{Def}[Thm]{Definition}
\newtheorem{Rem}[Thm]{Remark}

\newtheorem{Lem}[Thm]{Lemma}

\newtheorem{Con}[Thm]{Conjecture}

\def\M1{\ensuremath{{\mathcal M}_{1,1}}}

\def\CC{\ensuremath{\mathcal C}}

\def\EE{\ensuremath{\mathcal E}}
\def\FF{\ensuremath{\mathcal F}}

\def\MM{\ensuremath{\mathcal M}}

\def\OO{\ensuremath{\mathcal O}}

\def\TT{\ensuremath{\mathcal T}}

\begin{document}

\title{Semisimple Quantum Cohomology and Blow-ups}

\author{Arend Bayer}

\address{Arend Bayer,
Max-Planck-Institut f\"ur Mathmatik, Vivatsgasse 7, 53111 Bonn, Germany}

\email{bayer@mpim-bonn.mpg.de}

\begin{abstract}
Using results of Gathmann, we prove the following theorem: If a
smooth projective variety $X$ has generically semisimple $(p,p)$-quantum
cohomology, then the same is true for the blow-up of $X$ at any number
of points. This a successful test for a modified version of
Dubrovin's conjecture from the ICM 1998.
\end{abstract}

\maketitle


\section{Introduction}

This note is motivated by a conjecture proposed by Boris Dubrovin
in his talk at the ICM in Berlin 1998. It claims that the quantum cohomology
of a projective variety $X$ is generically semisimple if and only if its
bounded derived category $D^b(X)$ of coherent sheaves is generated by
an exceptional collection. We discuss here a modification of this conjecture
proposed in \cite{BM} and show its compatibility with blowing up at a
point.

Quantum multiplication gives (roughly speaking)
a commutative associative multiplication
$\circ_\omega \colon H^*(X) \otimes H^*(X) \to H^*(X)$ depending on
a parameter $\omega \in H^*(X)$. Semisimplicity of quantum cohomology means
that for generic parameters $\omega$, the resulting algebra is semisimple.
More precisely, quantum cohomology produces a formal Frobenius supermanifold
whose underlying manifold is the completion at the point zero of $H^*(X)$.
We call a Frobenius manifold generically semisimple if it is purely even
and the spectral cover
map $\Spec (\TT \MM, \circ) \to \MM$ is unramified over a general fibre.
Generically semisimple Frobenius manifolds are particularly well understood.
There exist two independent classifications of their germs, due to
Dubrovin and Manin. Both identify a germ via a finite set of invariants.
As mirror symmetry statements include an isomorphism
of Frobenius manifolds, this means that in the semisimple case one will have
to control only this finite set of invariants.

In \cite{BM}, it was proven that the even-dimensional part $H^{\text{ev}}$ of
quantum cohomology cannot be semisimple unless
$h^{p,q} = 0$ for all $p \neq q, p+q \equiv 0 \pmod 2$. On the other hand,
the subspace $\bigoplus_p H^{p,p}(X)$ gives rise to a Frobenius submanifold.
This suggested the following modification of Dubrovin's conjecture:
The Frobenius submanifold of $(p,p)$-cohomology is semisimple if and
only if there exists an exceptional collection of length
$\rk \bigoplus_p H^{p,p}(X)$.

A consequence of this modified conjecture is the following: 
\emph{Whenever $X$ has semisimple $(p,p)$-quantum cohomology,
the same is true for the blow-up of $X$ at any number of points.}
We prove this in Theorem \ref{semisimpleblowups}.

We would like to point out that our result suggests another small
change of the formulation of Dubrovin's conjecture. Dubrovin assumed that
being Fano is an additional necessary condition for semisimple quantum
cohomology.
However, as our result holds for the blow-up at an arbitrary number of
points, it yields many non-Fano counter-examples. We suggest to just drop any
reference to $X$ being Fano from the conjecture.

\subsection*{Acknowledgements}

I am grateful to Andreas Gathmann who made helpful comments on an early
version of this paper, in particular referring me to
J.~Hu's results.
I am indebted to my advisor Yuri~I.\ Manin for continued inspiring
and encouraging guidance and support.

\section{Definitions and Notations}

								\label{Defs}
							\label{FrobDefs}

Let $X$ be a smooth projective variety over
\C.  By $H_X := \bigoplus H^{p,p}(X, \C)$, we denote the
 space of $(p, p)$-cohomology.  Let $\Delta_0, \dots, \Delta_m,
\Delta_{m+1}, \dots, \Delta_r$ be a homogeneous basis of $H_X$,
such that $\Delta_0$ is the unit, and $\Delta_{m+1}, \dots, \Delta_r$
are a basis of $H^{1,1}(X)$.

We denote the correlators in the quantum cohomology of $X$ by
$$
\langle \Delta_{i_1} \dots \Delta_{i_n} \rangle_\beta.
$$
This is the number of appropriately counted stable maps
$$
f \colon (\CC, y_1, \dots, y_n) \to X
$$
where \CC{} is a semi-stable curve of genus zero,
$y_1, \dots, y_n$ are marked points on \CC,
the fundamental class of \CC{} is mapped to $\beta$ under $f$,
and $\Delta_{i_1}, \dots, \Delta_{i_n}$ are cohomology classes
representing conditions for the images of the marked points. In the case
of $\beta = 0$ it is artificially defined to be zero if $n \neq 3$, and
equal to $\int_X \Delta_{i_1} \cup \Delta_{i_1} \cup \Delta_{i_3}$ if
$n=3$.

Such a correlator vanishes unless
\begin{equation}					\label{dimaxiom}
k(\beta) := (c_1(X), \beta) = 3 - \dim X + \sum (\frac{|\Delta_{i_j}|}{2} - 1)
\end{equation}
where $|\Delta_{i_j}|$ are the degrees of the cohomology classes.

Before writing down the potential of quantum cohomology and the resulting
product, we will define the ring that it lives in. Let
$\{x_k | k \le m \}$ be the dual coordinates of $H_X/H^{1,1}(X)$
corresponding to the homogeneous basis $\{\Delta_k\}$. Instead of dual
coordinates in $H^{1,1}(X)$, we want to consider exponentiated coordinates.
This is done most elegantly by adjoining a formal coordinate
$q^\beta$ for effective classes $\beta \in H_2(X, \Z)/\text{torsion}$
with $q^{\beta_1 + \beta_2} = q^{\beta_1} q^{\beta_2}$.
Now let
$$
F_X = \Q[[x_k, q]]
$$
be the completion of the polynomial ring generated by $x_k$ and monomials
$q^\beta$ with  $\beta$ effective.

We consider $F_X$ as the structure ring of the formal Frobenius manifold
associated to $H_X$. The vector space $H_X$ acts on $F_X$ as a space
of derivations: $\Delta_k, k \le m$ acts as $\dbyd{x_k}$, and the
divisorial classes $\Delta_k, k > m$ act via
$q^\beta \mapsto (\Delta_k , \beta) q^\beta$. Hence we can formally consider
$H_X$ as the space of horizontal tangent fields of the formal Frobenius
manifold $\MM$, and $F_X \otimes H_X$ as its tangent bundle $\TT \MM$.

The flat structure of this formal manifold is given by the Poincar\'e pairing
$g$ on $H_X$. Given the flat metric, the whole structure of a formal Frobenius
manifold is an algebra structure on $F_X \otimes H_X$ over $F_X$ given
by the third partial derivatives of a potential $\Phi \in F_X$:
$$
g (\Delta_i \circ \Delta_j, \Delta_k) = \Delta_i \Delta_j \Delta_k \Phi
$$

To be able to consistently work only with exponentiated coordinates on
$H^{1,1}$, we slightly deviate from this definition: We use only
the non-classical part
$$
\Phi_X = \sum_{\beta \neq 0} \langle e^{\sum_{k \le m} x_k \Delta_k} \rangle_\beta
		q^\beta.
$$
of the Gromov-Witten potential (it is a consequence of the divisor axiom
that it makes sense to write $\Phi_X$ in this way), and define the product via
$g (\Delta_i \circ \Delta_j, \Delta_k) = g(\Delta_i \cup \Delta_j, \Delta_k)
+ \Delta_i \Delta_j \Delta_k \Phi_X$.
The choice of the ring $F_X$ is governed by the two properties that
it has to contain $\Phi_X$, and that $H_X$ has to act on it 
as a vector space of derivations. This is enough to ensure that all
standard constructions associated to a Frobenius manifold are defined over
$F_X$.

Explicitly, the multiplication is given by
\begin{equation}					\label{QCmult}
\Delta_i \circ \Delta_j = \Delta_i\,\cup\,\Delta_j +
\sum_{\beta \ne 0}\sum_{k\ne 0}
\langle\Delta_i\Delta_j\Delta_k e^{\sum_{k \le m} x_k \Delta_k}
\rangle_{\beta}\Delta^k q^\beta
\end{equation}
where $\Delta^k$ are the elements of the basis dual to $(\Delta_k)$ with
respect to the Poincar\'e pairing. The multiplication endows
$F_X \otimes H_X$ with the structure of a commutative, associative algebra
with $1 \otimes \Delta_0$ being the unit.

We call the whole structure of the formal Frobenius manifold on $F_X$ and
$H_X$ \emph{reduced quantum cohomology}.
The map of rings
\begin{equation} \label{Mspectrcover}
F_X \to F_X \otimes H_X, \quad f \mapsto f \otimes \Delta_0
\end{equation}
is the formal replacement of the spectral cover map 
$ \Spec (\TT \MM, \circ) \to \MM$
of a non-formal Frobenius manifold.

\begin{Def}
$X$ has semisimple reduced quantum cohomology if the
spectral cover map (\ref{Mspectrcover}) is generically unramified.
\end{Def}

More concretely, semisimplicity over a geometric point
$F_X \to k$ of $F_X$ means that after base change to $k$, the ring
$k \otimes H_X$ with the quantum product is isomorphic to
$k^{r+1}$ with component-wise
multiplication. Generic semisimplicity means that this is true for
a dense open subset in the set of $k$-valued points of $F_X$.

Finally, we recall the definition of the Euler field of quantum
cohomology. It is given by 
$$
\EE = -c_1(X)
+ \sum_{k \le m} \left(1 - \frac{|\Delta_k|}{2})\right) x_k \Delta_k.
$$

It induces a grading on $F_X$ and $F_X \otimes H_X$ by its Lie derivative.
E.\,g., a vector field is homogeneous of degree $d$ if
$\Lie_\EE (X) = [\EE, X] = d X$.
It is clear that the Poincar\'e pairing is homogenous of degree $(2 - \dim X)$
by the induced Lie derivative on $(H_X^*)^{\otimes 2}$. Further,
from the dimension axiom (\ref{dimaxiom}) it follows that $\Phi_X$ is
homogeneous of degree $(3 - \dim X)$. It a purely formal consequence of these
two facts that the multiplication $\circ$ is homogeneous of degree $1$ with
respect to $\EE$ (see \cite[I.2]{M}). 

\section{Semisimple quantum cohomology and blow-ups}

\subsection{Motivation}

So let us now assume that the variety $X$ satisfies the modified
version of Dubrovin's conjecture, i.\,e. that it has both an exceptional
collection of length $\rk \bigoplus_p H^{p,p}(X)$, and semisimple
reduced quantum cohomology. Let $\tilde X$ be its blow-up at some points.
By remark \ref{blowexmod}, this is a
test for the modified version of Dubrovin's conjecture \ref{Dubrpp}:
We know that $\tilde X$ has an exceptional system of desired length,
so it should have semisimple reduced quantum cohomology as well:

\begin{Thm}					\label{semisimpleblowups}
Let $\tilde X \to X$ be the blow-up of a smooth projective
variety $X$ at any number of closed points.

If the reduced quantum cohomology of $X$ is generically semisimple,
then the same is true for $\tilde X$.
\end{Thm}

In the case of dimension two, Del Pezzo surfaces were treated
in \cite{BM}, where the results of \cite{GP} on their quantum cohomology
were used.
The generalization presented here uses instead the results
in Andreas Gathmann's paper \cite{Gathmann}, with an improvement from the
later paper \cite{HuBlowup} by J.~Hu. The essential idea 
is a variant of the idea used in \cite{BM}: a partial
compactification of the spectral cover map where the exponentiated
coordinate of an exceptional class vanishes. However, in our case, this
is only possible after base change to a finite cover of the spectral
cover map.

\subsection{More notations}				\label{Morenots}

We want to compare the reduced quantum cohomology of $\tilde X$ with that
of $X$. We may and will restrict ourselves to the blow-up 
$j \colon \tilde X \to X$ of a single point. For the pull-back
$j^* \colon H^*(X) \to H^*(\tilde X)$ and the push-forward
$j_* \colon H^*(\tilde X) \to H^*(X)$ we have the identity
$j_* j^* = \id_{H^*(X)}$. Hence $H^*(\tilde X) = j^*(H^*(X)) \oplus \ker j_*$
canonically.  We will identify $j^* (H^*(X))$ with $H^*(X)$ from now on
and get a canonical decomposition 
$H_{\tilde X} = H_X \oplus H_E$ with
$ H_E =\bigoplus_{1\le k \le n-1} \C \cdot E^k$, where
$E$ is the exceptional divisor of $j$. The dual coordinates
$(x_k) =: \underline{x}$ on
$H_X/H^{1,1}(X)$ get extended via coordinates
$(x^E_2, \dots, x^E_{n-1}) =: \underline{x}^E$ to dual coordinates of
$H_{\tilde X}/H^{1,1}({\tilde X})$.
Let $E' \in H_2(\tilde X)$ be the class of a line in the exceptional
divisor $E \cong \P^{n-1}$. From Poincar\'e duality and the decomposition
of $H^*(\tilde X)$, we get a corresponding decomposition 
$H_2(\tilde X, \Z) = H_2(X, \Z) \oplus \Z \cdot E'$ in
homology, where we have identified $H_2(X)$ with its image via the dual
of $j_*$. With this identification, the cone
of effective curves in $X$ is a subcone of the effective cone in $\tilde X$.
Hence $F_X$ is a subring of $F_{\tilde X}$. We will call elements
$\beta \in H_2(X) \subset H_2(\tilde X)$ non-exceptional, and
$\beta \in \Z E'$ purely exceptional. 

We can also view $F_X$ as a quotient of $F_{\tilde X}$: 
Let $I$ be the completion of the subspace in $F_{\tilde X}$ generated by
monomials
$\underline{x}^{\underline{a}} \cdot (\underline{x}^E)^{\underline{b}}
q^{\tilde \beta}$ with $\underline{b} \neq (0, \dots, 0)$ or
$\tilde \beta \not \in H_2(X)$. Then evidently $F_X = F_{\tilde X}/I$. But
note that $I$ is \emph{not} an ideal, as there are effective classes
$\tilde \beta_1, \tilde \beta_2 \in H_2(\tilde X) \setminus H_2(X)$ whose sum
$\tilde \beta_1 + \tilde \beta_2$ is in $H_2(X)$.

Also, it is not true that $F_X \otimes H_X$
is a subring of $F_{\tilde X} \otimes H_{\tilde X}$. The next section
will summarize the results of \cite{Gathmann} that will enable us to study
the relation between the two reduced quantum cohomology rings.

\subsection{Gathmann's results}

\begin{Thm}
							\label{Gathmann}
The following assertions relate the Gromov-Witten invariants
of $\tilde X$ to those of $X$ (which we will denote by
$\langle \dots \rangle_\beta ^{\tilde X}$ and
$\langle \dots \rangle_\beta ^{X}$, respectively):

\begin{enumerate}
\item
\begin{enumerate}
\item 						\label{nonexcept}

Let $\beta \in H_2(\tilde X)$ be any non-exceptional homology class---so 
$\beta$ is any element of $H_2 (X)$---, and let
$T_1, \dots, T_m$ be any number of non-exceptional classes
in $H^*(X) \subset H^*(\tilde X)$, which we can identify with their preimages in
$H^*(X)$. Then it does not matter whether we compute the
following Gromov-Witten invariants with respect to $\tilde X$ or $X$:
$$
   \langle T_1 \otimes \dots \otimes T_m \rangle_\beta^{\tilde X}
   = \langle T_1 \otimes \dots \otimes T_m \rangle_\beta^{X}. 
$$

\item						\label{pureexcept}
Consider the Gromov-Witten invariants
$\langle T_1 \otimes \dots \otimes T_m \rangle_\beta^{\tilde X}$
with $\beta$ being purely
exceptional, i.\,e. $\beta = d \cdot E'$.

If any of the cohomology
classes $T_1, \dots, T_m$ are non-exceptional, the invariant is zero.
All invariants involving only exceptional cohomology classes can be
computed recursively from the following:
$$
   \langle E^{n-1} \otimes E^{n-1} \rangle_{E'}^{\tilde X}  = 1.
$$
They depend only on $n$.

\end{enumerate}

\item
\begin{enumerate}

\item 							\label{algo}
Using the associativity relations, it is possible to compute all
Gromov-Witten invariants of $\tilde X$ from those mentioned above
in~\ref{nonexcept} and~\ref{pureexcept}.

\item							\label{mixed}
\emph{Vanishing of mixed classes:}
Write $\tilde \beta \in H_2(\tilde X)$ in the form
$\tilde \beta= \beta + d \cdot E'$ where
$\beta$ is the non-exceptional part;
assume that $\beta \neq 0$. Let $T_1, \dots, T_m$ be non-exceptional
cohomology classes.
Let $l$ be a non-negative integer, and let
$2 \le k_1, \dots, k_l \le n-1$ be integers satisfying
\[ (k_1-1) + \dots + (k_l-1) < (d+1) (n-1). \]
Unless we have both $d = 0$ and $l = 0$, this implies the vanishing
of
\[ \langle T_1 \otimes \dots \otimes T_m \otimes
	   E^{k_1} \otimes \dots \otimes E^{k_l} \rangle_\beta = 0. \]

\end{enumerate}
\end{enumerate}
\end{Thm}

\begin{Prf}
The statement in no.~\ref{nonexcept} is proven by J.~Hu in
\cite[Theorem 1.2]{HuBlowup}. This is lemma 2.2 in \cite{Gathmann};
since the proof of this lemma is the only place where Gathmann uses
the convexity of $X$ (see remark 2.3 in that paper),
we can drop this assumption from his theorems.

The other equations follow trivially from statements in lemma 2.4
and proposition 3.1 in \cite{Gathmann}.
\end{Prf}

\subsection{Proof of Theorem \ref{semisimpleblowups}}

Let us first restate Gathmann's results in terms of the potentials
$\Phi_X$ and $\Phi_{\tilde X}$: We can write
$\Phi_{\tilde X}$ as
\begin{equation}			\label{potentials}
\Phi_{\tilde X} = \Phi_X + \Phi_{\text{pure}} + \Phi_{\text{mixed}}
\end{equation}
where $\Phi_X$ is the sum coming from all non-exceptional $\tilde \beta
=\beta$ and non-exceptional cohomology classes
(coinciding with the potential of $X$ by no.~\ref{nonexcept}),
$\Phi_{\text{pure}}$ is the sum coming from all correlators with
$\tilde \beta$ being purely exceptional (i.\,e. a positive multiple of $E'$), 
and $\Phi_{\text{mixed}}$ the sum from correlators
for mixed homology classes $\tilde \beta = \beta + d \cdot E'$ with
$0 \neq \beta \in H_2(X)$ and $d \neq 0$.


Now let $\tilde \EE$ and $\EE$ be the Euler fields of $\tilde X$ and $X$,
respectively. Let us consider the grading induced by
$$
\tilde \EE - \EE = (n-1) E + \sum_{2\le k \le n-1} (1 - k) x^E_k E^k.
$$

\begin{Lem}					\label{potdegrees}
With respect to $\tilde \EE - \EE$, the potential $\Phi_{\text{pure}}$
is homogeneous of degree $3 - n$, and
$\Phi_{\text{mixed}}$ only has summands of degree less than or equal to 
to $1-n$.
\end{Lem}

\begin{Prf}
The assertion about $\Phi_{\text{pure}}$ is just the dimension axiom
(\ref{dimaxiom}) of $\tilde X$, as $\EE \Phi_{\text{pure}} = 0$.
The statement about $\Phi_{\text{mixed}}$ is equivalent to
Gathmann's vanishing result, theorem~\ref{Gathmann} no.~\ref{mixed}.
\end{Prf}

Let $J \vartriangleleft F_{\tilde X}$ be the ideal generated by
$x_2^E, \dots, x_{n-1}^E$. We will show that the spectral cover map of
$\tilde X$ is already generically semisimple when restricted to the
fibre
\begin{equation}					\label{spectcover}
F_{\tilde X}/J \to H_{\tilde X} \otimes F_{\tilde X}/J.
\end{equation}

Write a monomial $q^{\tilde \beta}$ in $F_{\tilde X}$ as 
$q^{\tilde \beta} = Q^{-d}q^\beta$ if $\tilde \beta = \beta + d \cdot E'$
with $\beta \in H_2(X)$.
We make the base change to the cover given by adjoining
$Z:= \sqrt[n-1]{Q}$. More precisely, we first localize\footnote{Note that
$Q$ itself is not an element of $F_{\tilde X}$.} at $Q^{-1}$ and adjoin
a $(n-1)$-th root of $Q$:
We consider $R := \left(F_{\tilde X}/J\right)[Q][Z] / (Z^{n-1}\!-\!Q)$.

On the other hand, consider the subring $B$ of $R$ that consists of power
series in which $Z$ only appears with non-negative degrees.\footnote{
The ring $B$ is neither $F_X[[Z]]$ nor $F_X[Z]$; it is a different
completion of $F_X[Z]$.}
Then $R$ is a completion of the localiation $B[Z^{-1}]$ of $B$.
We claim that the quantum product ``extends'' to a product over $B$.

We define $M$ as the free $B$-submodule of
$B \otimes H^* (\tilde X)$ generated by
\[ \langle H^*(X), ZE, Z^2 E^2, \dots, Z^{n-1} E^{n-1} = Q E^{n-1}
				\rangle.			\]
More invariantly, $B$ is the completed subspace of $R$ generated by
monomials with non-positive degree with respect to $\tilde \EE - \EE$. And $M$
is the submodule of $B \otimes H_{\tilde X}$ generated
by $B \otimes H_X$ and all elements of strictly negative degree
in $B \otimes H_E$.

\begin{Lem}
\begin{itemize}
\item The quantum product restricts to $M$, i.\,e. 
$M \circ M \subseteq M$, and we have the following
cartesian digram:
\[
\xymatrix{
{B}
	\ar[d]
	\ar[r] &
{M}
	\ar[d]
			\\
{R}
	\ar[r] &
{R \otimes H^*(\tilde X)}
}	 \]

\item Now we take the push-out with respect to $B \to B/(Z) = F_X$. Then
the spectral cover map decomposes as 
\[ \xymatrix{
{B}
	\ar[r] \ar[d] &
{M}
	\ar[d] \\
{F_X} 
	\ar[r] &
{\left(F_X \otimes H_X\right)
 \oplus F_X[z]/(z^{n-1}- (-1)^{n-1})}
} \]
where the product on $F_X \otimes H_X$ is the quantum product of $X$.
\end{itemize}
\end{Lem}

First, we show how to derive Theorem~\ref{semisimpleblowups} from
the above lemma.  By the induction hypothesis, $F_X \to F_X \otimes H_X$
is generically semisimple. The second part of the lemma then tells us
that the map $B \to M$ is generically semisimple over the fibre of $Z=0$.

E.~g. by the criterion~\cite[IV, 17.3.6]{EGA} of unramifiedness,
it is clear that semisimplicity is an open condition for finite flat
maps. Hence, also $B \to M$ is generically semisimple.
The same is then true for its completed localization
$(F_{\tilde X}/J)[Q][Z]/(Z^{n-1}\!-\!Q)$. It is also evident that
the finite extension
$(F_{\tilde X}/J)[Q] \to (F_{\tilde X}/J)[Q][Z]/(Z^{n-1}\!-\!Q)$ cannot
change generic semisimplicity. Hence the spectral cover map
(\ref{spectcover}) must be generically semisimple (as its localization at
$Q$ is). And again by openness
of semisimplicity, it also holds for the full reduced quantum cohomology of
$\tilde X$.

\begin{Prf}[of the lemma]
We want to analyze the behaviour of multiplication with respect to the
grading of $\tilde \EE - \EE$. We decompose the quantum product
$\circ_{\tilde X}$, understood as a bilinear map
$(B \otimes H_{\tilde X}) \otimes (B \otimes H_{\tilde X}) \to
B \otimes H_{\tilde X}$,
into a sum
$\circ_{\tilde X} = \circ_{X} + \circ^E_{\text{class}} + \circ_{\text{pure}}
+ \circ_{\text{mixed}}$
according to the decomposition of $\Phi_{\tilde X}$ in (\ref{potentials});
we have written $\circ^E_{\text{class}}$ for the classcial cup product
of exceptional classes $E^i \circ^E_{\text{class}} E^j = E^{i+j}$ for
$0 \le i, j \le n-1$ and $i>0$ or $j>0$.
So for example $\circ_{\text{pure}}$ is defined by
$\tilde g (U \circ_{\text{pure}} V, W) = UVW \Phi_{\text{pure}}$ with
$\tilde g$ as the Poincar\'e pairing on $\tilde X$.

We claim that $\circ_X$, $\circ_{\text{pure}}$ and $\circ_{\text{mixed}}$ are
of degree $0$, $1$ and $\le -1$, respectively.

This is clear for $\circ_{X}$ and follows with standard Euler field
computations from the assertions in lemma \ref{potdegrees} (compare with
the computations in \cite[I.2]{M}):

Take a homogeneous component $\Phi_d$ of degree $d$ of any of the two
relevant potentials, and $\circ_d$ the corresponding component of the
multiplication. Let $U$, $V$ and $W$ be vector fields of degree
$u$, $v$ and $w$, respectively:
\begin{equation} \label{Euler1}
\begin{split}
(\tilde \EE - \EE) \tilde g(U \circ_d V, W)
&= (\tilde \EE - \EE) UVW \Phi_d
					\\
&= [\tilde \EE - \EE, U] V W \Phi_d + U [\tilde \EE - \EE, V] W \Phi_d
					\\
&\quad + UV [\tilde \EE - \EE, W] \Phi_d + UVW (\tilde \EE - \EE)\Phi_d
					\\
&= (u + v + w + d) UVW\Phi_d		\\
&= (u+v+w+d) \tilde g(U \circ_d V, W)
\end{split}
\end{equation}

Now write $\tilde g = g + g^E$ where $g$ is the Poincar\'e pairing of $X$
and $g^E$ the pairing of exceptional classes
$g^E(E^i, E^k) = \delta_{i+j,n} (-1)^{n-1}$. Then $g$ is of degree zero,
and $g^E$ of degree $2-n$ with respect to $\tilde \EE - \EE$. Let
$\circ_d = \circ_d^0 + \circ_d^E$ accordingly. Then
$U \circ_d V = U \circ_d^0 V + U \circ_d^E V$ is just the decomposition of
$U \circ_d V$ in the orthogonal sum
$H_{\tilde X} = H_X \oplus H_E$; in particular,
$U \circ_d V$ is homogeneous if and only if $U \circ_d^0 V$ and 
$U \circ_d^E V$ are.
So we have:
\begin{equation*} 
\begin{split}
(\tilde \EE - \EE) \tilde g(U \circ_d^0 V, W))
&= (\tilde \EE - \EE) g(U \circ_d^0 V, W) 				\\
&= g( [\tilde \EE - \EE, U \circ_d^0 V], W)
+ g(U \circ_d^0 V, [\tilde \EE - \EE,W]) 				\\
&= g( [\tilde \EE - \EE, U \circ_d^0 V], W) + w g(U \circ_d^0 V, W) \\
(\tilde \EE - \EE) \tilde g(U \circ_d^E V, W))
&= (\tilde \EE - \EE) g^E(U \circ_d^E V, W))			\\
&= \Lie_{\tilde \EE - \EE} (g^E)(U \circ_d^E V, W))		\\
&\quad + g( [\tilde \EE - \EE, U \circ_d^0 V], W)
  + g(U \circ_d^0 V, [\tilde \EE - \EE, W])			\\
&= g( (\tilde \EE - \EE) (U \circ_d^E V), W)
 + (2 - n + w) \tilde g(U \circ_d^E V, W)).
\end{split}
\end{equation*}
Comparing with (\ref{Euler1}), we see that 
$U \circ_d^0 V$ is of degree $u + v + d$, and $U \circ_d^E Y$ of degree
$u + v + d + n - 2$, in other words, $\circ_d^0$ has degree $d$ and
$\circ_d^E$ degree $d + n - 2$. Hence, the claim about the degree of
$\circ_{\text{mixed}}$ is obvious, and the one about
$\circ_{\text{pure}}$ follows from the additional fact the the derivative
of $\Phi_{\text{pure}}$ in $H_X$-direction is zero, so that
$\circ_{\text{pure}}^0$ is zero.

It is clear that $M$ is closed with respect to $\circ_X$ and
$\circ^E_{\text{class}}$.
That it is also closed under the multiplication
$\circ_{\text{mixed}}$ follows directly by degree reasons from the
description of $M$ in terms of degrees. With respect to
$\circ_{\text{pure}}$ we can argue via degrees if we additionally note that
$H_X\circ_{\text{pure}} H_{\tilde X} = 0$.

So we have proven $M \circ M \subseteq M$, and it remains to analyze the
product on $M/ ZM \cong F_X \otimes H_{\tilde X}$. Note that all elements
in $M$ of degree $\le -2$ are mapped to zero in this quotient.

It is clear that $\circ_X$ induces the quantum product of $X$ on the
subspace $F_X \otimes H_X$ and is zero on $H_E$. We already noted that
$H_X \circ_{\text{pure}} H_{\tilde X} = 0$.
Also, $Y_1 \circ_{\text{mixed}} Y_2$ is always
zero if $Y_1$ or $Y_2$ is in $M \cap B \otimes H_E$ for degree reasons.

We investigate the product with $ZE$. For this we can ignore $\circ_X$
and $\circ_{\text{mixed}}$. The classical part contributes
$ZE \circ^E_{\text{class}} (ZE)^{i} = (ZE)^{i+1}$ for $0 \le i \le n-1$.
For $\circ_{\text{pure}}$ we finally have to use the explicit
multiplication formula:
$$
ZE \circ_{\text{pure}} (ZE)^{i}
= (-1)^{n-1} \sum_{d > 0} \sum_j \langle E E^i E^j \rangle_{d E'}
			Z^{i+1} E^{n-j} Q^{-d}.
$$
By the dimension axiom, this can only be non-zero if
$(n-1) d = 3 - n + (1 - 1) + (i - 1) + (j-1)$, or, equivalently,
$(n-1)(d+1) = i + j$. This is only possible for $d=1$ and $i = j = n-1$,
where we have
$\langle E E^{n-1} E^{n-1} \rangle_{E'}
= -\langle E^{n-1} E^{n-1} \rangle_{E'} = -1$. 
We thus get
$$
ZE \circ (ZE)^{i} = 
\begin{cases}
       Z^{i+1}E^{i+1} & \text{if $i \le n-2$} \\
       (-1)^n Z E     & \text{if $i= n-1$.}
 \end{cases}
$$

Let $Y:= (-1)^n Q E^{n-1} = (-1)^n Z^{n-1} E^{n-1}$. 
As a consequence of the last equation, multiplication by $Y$ in the ring
$M / ZM$ is the identity on $(M \cap B \otimes H_E) / ZM \cong F_X \otimes H_E$.
In particular, $Y$ is an idempotent and gives a splitting of
$M/ZM \cong F_X \otimes H_E  \oplus K$ into the image
$F_X \otimes H_E$ and the kernel $K$ of $Y \circ$. The algebra structure
on $F_X \otimes H_E$ is isomorphic to $F_X[z]/(z^{n-1}-(-1)^{n-1})$ via
$z \mapsto ZE$.

The kernel is generated by $\Delta_1, \dots, \Delta_m, \Delta_0 - Y$,
and $\Delta_0 - Y$ is its unit. Mapping each element in $K$ to its degree
zero component, we get an isomorphism $K \to F_X \otimes H_X$ that maps
the multiplication on $K$ isomorphically to its degree zero component
$\circ_X$, and the lemma is proven.
\end{Prf}

\subsection{Further Questions}

The first example where our theorem applies is the case of $X = \P^n$.
For $n=2$, this yields the semisimplicity of quantum cohomology for all
Del Pezzo surfaces as proven earlier in \cite{BM}. Further, semisimplicity
has been established in \cite{tianxu-ss} by Tian and Xu, using results of
Beauville (see \cite{Beauville-SS}), for low degree complete intersections
in $\P^n$.

Generally speaking, once the three-point Gromov-Witten correlators are
known, and thus generators and relations for the small quantum cohomology
ring, it is an excercise purely in commutative algebra to check generic
semisimplicity in small quantum cohomology. For example, using Batyrev's
formula for Fano toric varieties \cite{BatyrevQCtoric} and its explicit
version for the projectivization of spliting bundles over $\P^n$ given in
\cite{AnconaMaggesiBundles}, semisimplicitiy can be shown to hold for
these bundles.

Of course, our theorem~\ref{semisimpleblowups} covers only the first
part of Dubrovin's conjecture. It would be very encouraging if it
was possible to show his statement on Stokes matrices in a similar
way. To my knowledge, the only case where this part has been
checked is the case of projective spaces (cf. \cite{Guzz}).

Revisiting Gathmann's algorithm to compute the invariants
of $\tilde X$ (Theorem~\ref{Gathmann},
no.~\ref{algo}), we notice that all the initial data it uses is already
determined by the multiplication in the special fibre $Z=0$
of our partially compactified spectral cover map. In other words,
the Frobenius manifold on $F_{\tilde X}$ and $H_{\tilde X}$ 
is already determined by the structure at $Z=0$.

Yet our construction does not yield a Frobenius structure at the 
divisor $Z=0$. If there was a formalism of Frobenius manifolds
with singularities along divisors, and if there was a
way to extend Dubrovin's Stokes matrices to these divisorial Frobenius
manifolds, this might also lead to an elegant treatment of Stokes
matrices of blow-ups.

Also, one would like to extend the method to the case of the blow-up along
a subvariety, analogously to Orlov's theorem~\ref{blowex}. The
next-trivial case of the blow-up along a fibre $\{x_0\} \times Y$ in a
product $X \times Y$ follows from our result and the discussion of products
in section \ref{DubrAProds}.

\section{Exceptional systems and Dubrovin's conjecture}
							\label{exceptsys}

In this section, we briefly review Dubrovin's conjecture and its modified
version, and explain how our theorem fits into
this context.

\subsection{Exceptional systems in triangulated categories}

We consider a triangulated category $\CC$. We assume that it is linear over
a ground field $\C$.

\begin{Def}				\label{Defexceptsys}
\begin{itemize}
\item An exceptional object in \CC{} is an object $\EE$  such that
the endomorphism complex of $\EE$ is concentrated in degree zero and
equal to $\C$:
\[ \RHom^\bullet (\EE,\EE) = \C [0] \]
\item An exceptional collection is a sequence $\EE_0, \dots, \EE_m$ of
exceptional objects, such that for all $i>j$ we have no morphisms
from $\EE_i$ to $\EE_j$:
\[ \RHom^\bullet (\EE_i, \EE_j) = 0 \quad \text{if $i>j$} \]
\item An exceptional collection of objects is called a complete
exceptional collection (or exceptional system), if the objects
$\EE_0, \dots, \EE_m$ generate \CC{} as a triangulated category:
The smallest subcategory of \CC{} that contains all $\EE_i$, and is closed
under isomorphisms, shifts and cones, is \CC{} itself.
\end{itemize}
\end{Def}

The first example is the bounded derived category $D^b(\P^n)$ on a projective
space with the series of sheaves $\OO(i), \OO(i+1), \dots, \OO(i+n)$ (for
any $i$). Exceptional systems were studied extensively by a group
at the Moscow University,
see e.~g. the collection of papers in \cite{Rudakov}.

More generally, exceptional systems exist on flag varieties; other
examples include quadrics in $\P^n$ and projective bundles over a variety
for which the existence of an exceptional system is already known.

\subsection{Dubrovin's conjecture}

On the other side of Dubrovin's conjecture we consider the Frobenius
manifold $M$ associated (as in \cite{M} or \cite{DubrPain}) to the quantum
cohomology of $X$.
As already mentioned in the introduction, Dubrovin's conjectures relates
generic semisimplicity of $M$ to the existence of an exceptional
system: 

\begin{Con}\cite{DubrICM}		\label{DubrConj}
Let $X$ be a projective variety.

The quantum cohomology of $X$ is generically semisimple if and only
if there exists an exceptional system in its derived category
$D^b(X)$.
\end{Con}

In further claims of his conjecture, he relates invariants of $M$ to
characteristics of the exceptional system: The so-called \emph{Stokes
matrix} $S$ of the Frobenius manifold should have entries
$S_{ij} = \chi (\EE_i, \EE_j)$. We almost completely omit these parts of his
conjecture in our discussion.

An expectation underlying Dubrovin's conjecture is that the
mirror partner of such a variety $X$ will be the unfolding of a
function with isolated singularities. The quantum cohomology should be
isomorphic to a Frobenius manifold structure on the base space of
the unfolding, as established by Barannikov for projective spaces,
cf.~\cite{BarPn}.

If $X$ has cohomology with Hodge indices other than $(p, p)$, it can neither
have an exceptional system, nor can the Frobenius manifold of its
quantum cohomology be semisimple:
\begin{itemize}
\item 
To make sense of all parts of Dubrovin's conjecture, an exceptional collection
should have length $\rk H^{\text{ev}}(X)$. But the length of
an exceptional collection is bounded by the rank of $N^*(X)$, the group
of algebraic cycles modulo numerical equivalence.\footnote{From the
Hirzebruch-Riemann-Roch theorem, it follows easily that the Chern
characters of the
exceptional objects are linearly independent.} And we always have
$\rk N^*(X) \le \rk H_X$.
\item
The subspace $\bigoplus_p H^{p,p}(X) \subset H^*(X)$ gives
rise to a Frobenius submanifold $M'$ of $M$; this is the Frobenius manifold
we constructed in section \ref{Defs}. This is
a maximal Frobenius submanifold of $M$ that has a chance of being
semisimple (\cite[1.8.1]{BM}).
\end{itemize}

This suggested the following modification:
\begin{Con}\cite{BM} \label{Dubrpp}
The variety $X$ has generically semisimple reduced quantum cohomology (i.\,e.
$M'$ is generically semisimple) if and only if there
exists an exceptional collection of length $\rk \bigoplus_p H^{p,p}(X)$ in
$D^b(X)$.
\end{Con}

\subsection{Products}				\label{DubrAProds}

It follows easily from well-known facts that Dubrovin's conjecture is
compatible with products, i.\,e. when it is true for two varieties $X, Y$,
it will also hold for their product $X \times Y$.

\begin{Thm}
Let $\EE_0, \dots, \EE_m$ be an exceptional system on $X$, and
$\FF_0, \dots, \FF_{m'}$ one on $Y$. Then $(\EE_{i_k} \boxtimes \FF_{j_k})_k$
forms an exceptional system on $X \times Y$,
where $(i_k, j_k)_k$ indexes the set
$\{1, \dots, m\} \times \{1, \dots, m'\}$ in any order such that we never
have $i_k > i_{k'}$ and $j_k > j_{k'}$ for $k < k'$.
\end{Thm}

This follows from the Leray spectral sequence computing the Ext-groups
on $X \times Y$. It also shows that the Stokes matrix of the exceptional
system on $X \times Y$ is the tensor product of the Stokes matrices
on $X$ and $Y$:
\[ \chi (\EE_{i_k} \boxtimes \FF_{j_k}, \EE_{i_{k'}} \boxtimes \FF_{j_{k'}})
 = \chi (\EE_{i_k}, \EE_{i_{k'}}) \cdot \chi (\FF_{j_k}, \FF_{j_{k'}})
\]

The corresponding statements hold for quantum cohomology:
Let $M$ and $M'$ be the Frobenius manifolds associated to the quantum
cohomology of $X$ and $Y$, respectively.
The Frobenius manifold of the quantum cohomology of $X \times Y$ is the
tensor product $M \otimes M'$ (\cite{KoMaProd}, \cite{BehrendProd},
\cite{KaKuenn}).
A pair of semisimple points in $M$ and $M'$ yields a semisimple point
in $M \otimes M'$, and the Stokes matrix of the tensor product is the
tensor product of the Stokes matrices of $M$ and $M'$ (\cite[Lemma
4.10]{DubrPain}).
It is also clear that the same holds for the reduced quantum cohomology
on $H_X$, $H_Y$ and $H_X \otimes H_Y$.

Hence, Dubrovin's conjecture follows for the product 
if it holds for $X \times Y$.
And in cases where
$H_X \otimes H_Y= H_{X\times Y}$, i.\,e.
$\bigoplus_p H^{p,p} (X) \otimes \bigoplus_p H^{p,p} (Y) = 
\bigoplus_p H^{p,p} (X \times Y)$,
the same holds for the modified conjecture \ref{Dubrpp}.

\subsection{Complete exceptional systems and blow-ups}
\label{exblowup}

\begin{Thm}\cite{Orlovblowup}				\label{blowex}
Let $Y$ be a smooth subvariety of the smooth projective variety
$X$. Let $\rho \colon \tilde X \to X$ be the blow-up of $X$ along $Y$.

If both $Y$ and $X$ have an exceptional system, then the same is true
for $\tilde X$.
\end{Thm}

Consider the case where $Y$ is a point:
Let $\P^{n-1}\cong E \subset \tilde X$ be the exceptional divisor
($n$ is the dimension of~$X$). If $\EE_0, \dots, \EE_r$ is a given
exceptional system in $D^b(X)$, then 
$\OO_E(-n+1), \dots, \OO_E (-2), \OO_E(-1), \rho^* (\EE_0), \dots,
\rho^*(\EE_r)$ is an exceptional system in $D^b(\tilde X)$. Hence, the
following holds:
\begin{Rem}						\label{blowexmod}
If $X$ has an exceptional collection of length $\rk H_X$,
then the analogous statement is true for the blow-up of $X$ at any number
of points.
\end{Rem}

\addcontentsline{toc}{section}{References}
\bibliography{all}                      
\bibliographystyle{alphaspecial}     

\end{document}